\documentclass[11pt,oneside]{amsart}

\usepackage{amscd}
\usepackage{amssymb}
\usepackage[latin1]{inputenc}
\usepackage[all]{xypic}
\usepackage{verbatim}

\def\Z{\mathbb{Z}}
\def\1{^{-1}}

\def\G{\Gamma}

\newtheorem{thm}{Theorem}

\title{On extensions of group with infinite conjugacy classes, I}

\begin{document}
\maketitle
\begin{center}
{\sc Jean-Philippe PR\' EAUX}\footnote[1]{Centre de Recherche de l'Armée de l'air, Ecole de l'air, F-13661 Salon de
Provence air}\ \footnote[2]{Centre de Math\'ematiques et d'Informatique, Universit\'e de Provence, 39 rue
F.Joliot-Curie, F-13453 marseille
cedex 13\\
\indent {\it E-mail :} \ preaux@cmi.univ-mrs.fr}
\end{center}

\begin{abstract}
We give a characterization of the group property of being with infinite conjugacy classes (or {\it icc}, {\it i.e.}
$\not= 1$ and of which all
conjugacy classes beside 1 are  infinite) for   extensions of abelian, centerless, icc, or word hyperbolic  groups.\\
\end{abstract}

\section*{Introduction}

A group is said to be with {\sl infinite conjugacy classes} (or {\it icc}) if it is non trivial, and if all its
conjugacy classes beside $\{ 1\}$ are infinite. This property is motivated by the theory of Von Neumann algebra, since
for any group $\G$, a necessary and sufficient condition for its Von Neumann algebra $W_\lambda^*(\G)$ to be a type
$II-1$ factor is that $\G$ be icc (cf. \cite{roiv}).

The property of being icc has been characterized in several classes of groups : 3-manifolds and $PD(3)$ groups in
\cite{aogf3v}, groups acting on Bass-Serre trees in \cite{ydc}, wreath products, semi-direct products and finite
extensions in \cite{p1,p2,p3}. We will focus here on groups defined by some specific extensions of a  group, namely
extensions of abelian groups, of centerless groups, of icc groups, or of word hyperbolic groups. We give a necessary
and sufficient condition for those groups to satisfy the icc property.\smallskip

We first consider the case of extensions of abelian groups. In a group $G$, $FC(G)$ denotes the unions of finite
conjugacy classes in $G$ ; it turns out that $FC(G)$ is a characteristic subgroup of $G$. Obviously, $G$ is icc if and
only if $FC(G)=1$.
\begin{thm}\label{t1} Let $G\not=1$ be a group defined by an extension of an abelian group :
$$1\longrightarrow K \text{abelian}\longrightarrow G\longrightarrow Q \longrightarrow 1$$
and $\theta : Q\longrightarrow Aut(K)$ the associated homomorphism. Then $G$ is icc if and only if both :
\begin{itemize}
\item[$(i)$] $K\setminus 1$ contains only infinite $\theta(Q)$-orbits, \item[$(ii)$] the induced homomorphism $\theta :
FC(Q)\longrightarrow Aut(K)$ is injective.
\end{itemize}
\end{thm}

On the opposite, we consider the case of extensions of centerless groups. For any extension $G$ of a group $K$ we
consider the homomorphism $\pi:G\longrightarrow Aut(K)$ defined by,  for all $g\in G$, $k\in K$, $\pi(g)(k)=k^g$.
\begin{thm}\label{t2}
Let $G\not=1$ be a group defined by an extension of a centerless group :
$$1\longrightarrow K\text{centerless}\longrightarrow G\longrightarrow Q\longrightarrow 1$$
and $\theta : Q\longrightarrow Out(K)$ the associated homomorphism. Then $G$ is icc if and only if both :
\begin{itemize}
\item[$(i)$] $K$ does not contain a non trivial normal subgroup $N$ either finite or $\Z^n$ which is preserved by
$\pi(G)$ and contains  only finite $\pi(G)$-orbits, \item[$(ii)$] the induced homomorphism $\theta :
FC(Q)\longrightarrow Out(K)$ is injective.
\end{itemize}

\end{thm}

Since in condition $(i)$ $N$ necessarily lies in $FC(K)$, and since an icc group is centerless, it follows immediatly :
\begin{thm}
Let $G$ be a group defined by an extension of an icc group :
$$1\longrightarrow K\text{icc}\longrightarrow G\longrightarrow Q\longrightarrow 1$$
and $\theta : Q\longrightarrow Out(K)$ the associated homomorphism. Then $G$ is icc if and only if the induced
homomorphism $\theta :FC(Q)\longrightarrow Out(K)$ is injective.
\end{thm}

Together with a few observations on word hyperbolic groups, one obtains as another corollary the following
characterization of icc property for extensions of word hyperbolic groups.
\begin{thm}\label{t3}
Let $G$ be defined by an extension of a non trivial word hyperbolic group :
$$1\longrightarrow K \text{hyperbolic}\longrightarrow G\longrightarrow Q\longrightarrow 1$$
and $\theta : Q\longrightarrow Out(K)$ the associated homomorphism. Then $G$ is icc if and only if both :
\begin{itemize}
\item[$(i)$] $K$ is icc, \item[$(ii)$] the induced homomorphism $\theta : FC(Q)\longrightarrow Out(K)$ is injective.
\end{itemize}
\end{thm}

\begin{comment}
\begin{thm}\label{t3}
Let $G$ be  defined by an extension of group :
$$1\longrightarrow K\longrightarrow G\longrightarrow Q\longrightarrow 1$$
and $\theta : Q\longrightarrow Out(K)$ the associated homomorphism. It gives rise to the homomorphism
$\pi:G\longrightarrow Aut(K)$ as the natural extension of $Inn(K)$ by $\theta(Q)$ ; suppose moreover that this last
extension splits.  Then $G$ is icc if and only if :
\begin{itemize}
\item[$(i)$] $FC(K)\setminus 1$ contains only infinite $\theta(Q)$-orbits, \item[$(ii)$] the induced homomorphism
$\theta : FC(Q)\longrightarrow Out(K)$ is injective.
\end{itemize}
\end{thm}
\end{comment}
\section{Proof of the results.}

Let $G$ be a group, $H$ a non empty subset of $G$, and $u,g\in G$. We use the notations $u^g=g^{-1}ug$ and
$u^H=\{u^g|g\in H\}$, so that $g^G$ denotes the conjugacy class of $g$ in $G$ ; $Z_G(g)$ stands for the centralizer of
$g$ in $G$ and $Z(G)$ for the center of $G$. It follows immediately that $g^G$ is finite if and only if $Z_G(g)$ has a
finite index in $G$. Let $\rho : G\longrightarrow Q$ denote the projection given by the extension.
\medskip

\noindent{\sl Proof of theorem \ref{t1}.} Note that that the conjugacy class in $G$ of $k\in K$, coincide with its
$\theta(Q)$-orbit. Obviously $G$ icc implies that both conditions $(i)$ and $(ii)$ are satisfied. We prove the converse
by showing that if $G$ is not icc then either condition $(i)$ or condition $(ii)$ is not satisfied. So suppose that $G$
is not icc ; since $G\not=1$ there exists $u\not=1$ in $G$ such that its conjugacy class $u^G$ in $G$ is finite.
If $u\in K$ ; then condition $(i)$ is not satisfied, so suppose in the following that $u\in G\setminus K$. Then
$\rho(u)$ lies in $FC(Q)$, for $\rho(u)^Q=\rho(u^G)$, and $\rho(u)\not=1$. Let $k\in K$, then $[u,k]\in K$ has a finite
conjugacy class in $G$, for $Z_G(u)$ and $Z_G(kuk^{-1})$ both have a finite index in $G$ and their intersection lies in
$Z_G([u,k])$. So that either, for all $k\in K$, $[u,k]=1$, which implies that $\pi(u)=\theta(\rho(u))=Id_K$ and then
$\theta :FC(Q)\longrightarrow Aut(K)$ is non injective, or as above, condition $(i)$ is not
satisfied.\hfill$\square$\medskip
%\in Z_G(u)\cap K$, then $u=g^{-1}ug=g^{-1}\theta(g)u$ so that $\theta(g)=g$.

\noindent{\sl Proof of theorem \ref{t2}.} Suppose that $G$ is not icc ; since $G\not=1$ there exists $u\not=1$ in $G$
such that $u^G$ is finite. If $u$ lies in $K$,  let $N'$ be the subgroup of $K$ generated by $u^G$ ; $N'$ is a non
trivial subgroup of $FC(K)$ preserved by $\pi(G)$ and with only finite $\pi(G)$-orbits. Moreover $N'$ is a finitely
generated $FC$-group and hence  a finite central extension of a f.g. abelian group (cf. \cite{fc}) ; it follows that
the set $Tor(N')$ of torsion elements of $N'$ is a characteristic finite subgroup of $N'$ and $N'/Tor(N')$ is $\Z^n$.
Hence condition $(i)$ is not satisfied by picking $N=Tor(N')$, if $Tor(N')\not=1$, and $N=N'$ otherwise.
 If $u$ lies
in $G\setminus K$, then $\rho(u)$ is a non trivial element of $FC(Q)$. For any $k\in Q$, $[u,k]$ lies in $K$ and has a
finite conjugacy class in $G$, so that, as in the proof of theorem \ref{t1}, we can suppose that $\pi(u)$ is the
identity, and hence $\theta(\rho(u))=1$, for otherwise condition $(i)$ is not satisfied. It follows that
$\theta:FC(Q)\longrightarrow Out(K)$ is non injective, so that condition $(ii)$ is not satisfied. We have proved that
if both condition $(i)$ and $(ii)$ are satisfied then $G$ is icc ; we now prove the converse. If $K$ contains a non
trivial subgroup preserved by $\pi(G)$ and with only finite $\pi(G)$-orbits then each element of $N$ has a finite
conjugacy class in $G$ and then $G$ is not icc.  If $\theta:FC(Q)\longrightarrow Out(K)$ is non injective, let
$q\not=1$ in $FC(Q)$ be such that $\theta(q)=1$. Then there exists a lift $\bar{q}$ of $q$ in $G$ such that
$\pi(\bar{q})$ is the identity of $K$, so that $K\subset Z_G(\bar{q})$. For any $q_1\in Z_Q(q)$, if $\bar{q_1}$ denotes
a lift of $q_1$ in $G$, then $\bar{q_1}\bar{q}\bar{q_1}^{-1}=\bar{q}k$ for some $k\in K$, which must belong to $Z(K)$
since $\pi(\bar{q_1}\bar{q}\bar{q_1}^{-1})=\pi(k)$ is the identity on $K$. Since $Z(K)=1$, $Z_G(\bar{q})$ contains
$\rho^{-1}(Z_Q(q))$, and hence has a finite index in $G$, so that $G$ is not icc. \hfill$\square$\medskip

\noindent{\sl Proof of theorem \ref{t3}.} If $K$ is a non trivial elementary group then $K$ is not icc, and neither
$G$, for $Aut(K)$ is finite. So that the assumption in theorem \ref{t3} is satisfied and we suppose in the following
that $K$ is non elementary.
 A non elementary word hyperbolic has a finite center (cf. \cite{gro}). It follows that
each element of $Z(K)$ has a finite conjugacy class both in $K$ and in $G$. Then, if $Z(K)\not=1$ neither $K$ nor $G$
are icc so that the assumption of theorem \ref{t3} is satisfied. So suppose in the following that $K$ is centerless, so
that theorem \ref{t2} applies. In a non elementary word hyperbolic group $K$, $FC(K)$ turns out to be a finite
characteristic subgroup of $K$. Hence theorem \ref{t2}  shows that $G$ is icc if and only both
$\theta:FC(Q)\longrightarrow Out(K)$ is injective and $FC(K)=1$.\hfill$\square$

\end{document}